\documentclass{amsart}
\usepackage{latexsym}
\usepackage{amsfonts}
\usepackage{amsmath}
\usepackage{amssymb}
\usepackage{tikz} 
\usetikzlibrary{shapes}
\usetikzlibrary{shapes.geometric}
\usepackage{color}
\usepackage{pict2e}

\newtheorem{theorem}{Theorem}[section]
\newtheorem{lemma}[theorem]{Lemma}

\theoremstyle{definition}

\theoremstyle{remark}

\numberwithin{equation}{section}

\newcommand{\R}{\mathbb{R}}

\newcommand{\clx}{CL({X})}
\newcommand{\cx}{C({X})}

\newcommand{\cnx}{C_{n}(X)}

\newcommand{\cpx}[1]{C_{p}(X_{#1})}

\newcommand{\Hx}{\mathcal{H}(X)}

\newcommand{\graphs}{\mathcal{X}}
\newcommand{\homeo}{\approx}
\newcommand{\haus}{\tau_{H}}
\newcommand{\vietoris}{\tau_{V}}

\begin{document}


\vspace{0.5in}

\title{The hyperspaces $\cnx$ for finite ray-graphs}

\author{Norah Esty}
\email{nesty@stonehill.edu}
\address{Department of Mathematics, Stonehill College,
Easton, Massachusetts 02357}



\begin{abstract} 
In this paper we consider the hyperspace $\cnx$ of non-empty and closed subsets of a base space $X$ with up to $n$ connected components.  We consider a class of base spaces called finite ray-graphs, which are a noncompact variation on finite graphs.  We prove two results about the structure of these hyperspaces under different topologies (Hausdorff and Vietoris). \\
\end{abstract}

\maketitle

\section{Introduction}

The last thirty years has produced a large amount of research in the area of hyperspaces.  A hyperspace is a topological space whose points are subsets of a given base space.  A general hyperspace is denoted $\mathcal{H}(X)$, where $X$ is the base space.  Common hyperspaces include $CL(X)$, the space of non-empty and closed subsets of $X$, $C_{n}(X)$, the space of non-empty, closed subsets of $X$ with up to $n$ connected components (where $C(X) = C_{1}(X)$), and $F_{n}(X)$, the space of non-empty subsets with up to $n$ points (called the ``$n$th symmetric product."). There are several topologies available for such spaces.  For $CL(X)$, if the base space is compact, two of the most popular topologies, the Hausdorff and the Vietoris, agree.   However, when the base space is not compact, they differ, and in fact the Vietoris topology is non-metrizable.  In contrast, by using a bounded metric on the space, or allowing for infinite-valued metrics, the Hausdorff topology arises from a metric.  Most of the study of hyperspaces has been done in the case where the base space $X$ is a continuum. \\

In 1968, Duda did an examination of the hyperspace of subcontinua of finite connected graphs, and under some minor conditions was able to give a description of $C(X)$ as a polyhedron, decomposable into balls of various dimensions.  \cite{Duda1}, \cite{Duda2}.  A single hyperspace may consist of several sections of different dimension: a two-dimensional disc glued to a three dimensional ball, etc.    In particular, for $X$ a finite graph, the hyperspace $C(X)$ is known to be compact and connected. \\

Uniqueness of hyperspaces is the property that if $\mathcal{H}(X)$ is homeomorphic to $\mathcal{H}(Y)$, then $X$ is homeomorphic to $Y$.  This is not true in general, so the question has become for which classes of base spaces it holds.  Work by Acosta, Duda, Eberhart, and Nadler has shown that finite graphs (different from an arc and the simple closed curve), hereditarily indecomposable continua, and smooth fans have unique hyperspace $\cx$.  See \cite{Duda1}, \cite{Duda2}, \cite{Acosta1}, \cite{Acosta2}, \cite{EberhartNadler}.  In 2002 and 2003, Illanes continued this study, and showed that for finite graphs the hyperspaces $C_{n}(X)$ are unique.  See \cite{Illanes1}, \cite{Illanes2}.  \\

In this paper we are interested in the situation where the base space is not compact.  We look at a natural generalization of finite graphs which we call finite ray-graphs, which consist of vertices, edges, and rays.  Because the graphs are not compact we must always specify which topology we are using, and in section~\ref{sec:hausdorff} of this paper we will use first the Hausdorff topology (arising from the Hausdorff metric, which we allow to be infinite-valued) and in section~\ref{sec:vietoris} the Vietoris.  \\

To assist the reader, in sections~\ref{sec:examples} and \ref{sec:joins} we present many models of the hyperspace $C(X)$.  In section~\ref{sec:joins} we state a theorem about hyperspace $C(X \vee_{p} Y)$ of a wedge product at a point, when the hyperspaces of the $C(X)$ and $C(Y)$ are known.  We state this theorem without proof, as it seems to be well-known in the folk-lore (although we have been unable to find a reference).  This theorem gives a nice algorithm for drawing hyperspaces.  \\

In sections~\ref{sec:hausdorff} and~\ref{sec:vietoris} we prove two main results about the number of connected components of the hyperspace $C_{n}(X)$ of a finite, connected ray-graph $X$: once in the Hausdorff topology and once in the Vietoris.   In particular, we show that when allowing for an infinite-valued Hausdorff metric, a finite, connected ray-graph with $k$ rays will have a hyperspace $C_{n}(X)$ with $2^{k}$ connected components for all $n$, and will not be compact.   In contrast, under the Vietoris topology $C_{n}(X)$ is connected for all $n$.

\section{Preliminaries and Notation}
\label{sec:intro}
\subsection{Notation}
\label{sec:notation}

There is not always consistent notation used for the different hyperspaces of a given base space $X$.  We attempt to use those notations from the literature which are least ambiguous.  Given a metric space $X$, we define the following notation for the hyperspaces we will discuss:

\begin{itemize}
\item $\clx = \{A \subset X: A$ closed and $A \neq \emptyset \}$
\item $\cnx = \{A \subset X: A$ closed, $A \neq \emptyset$ and $A$ has at most $n$ connected components $\}$
\item $\cx = C_{1}(X)$
\end{itemize}

It should be pointed out that much of the literature on hyperspaces assumes that the base space $X$ is compact, in effect making $\cx$ the hyperspace of subcontinua, but we are not assuming that here.  This is also why we write $\clx$ rather than $2^{X}$, which is more common, but to many readers may mean \emph{bounded} closed subsets, which we do not mean.   When we wish to refer to a general hyperspace,  we will write $\Hx$.  \\

Initially we will endow our hyperspaces with the Hausdorff topology ($\haus$).  The Hausdorff topology has the virtue that is arises from a metric, although since we are interested in unbounded base spaces, we allow the metric to be infinite-valued.  Let $\mathcal{H}(Y)$ be a hyperspace over a metric base space $Y$.  If $A, B \subset Y$, and if $N_{Y}(A, \epsilon)$ indicates the $\epsilon$-neighborhood in the space $Y$ around the subset $A$,  then the Hausdorff distance in the hyperspace is given by 

$$d_{H}(A, B) = \inf\{\epsilon: A \subset N_{Y}(B, \epsilon) \mbox{ and } B \subset N_{Y}(A, \epsilon) \}$$

If the elements of the hyperspace are not closed subsets, then it is possible to have the distance between two non-equal sets be zero.  However we will deal exclusively with closed sets.  One can see from this definition that if $A$ is bounded and $B$ is not, the Hausdorff distance between $A$ and $B$ is infinite. \\

Later in the paper we will use the Vietoris topology.  This topology is usually given by a basis or subbasis definition, which we will recall in section~\ref{sec:vietoris}.

\subsection{The class of base spaces: finite ray-graphs}
\label{sec:basespaces}

For our base spaces, we will consider a variation on finite graphs, which we will call finite ray-graphs.  These graphs will consist of a finite number of vertices (points), edges (homeomorphic to $[0,1]$ and attached at two vertices, or at one vertex twice) and rays (homeomorphic to $[0, \infty)$ and attached at one vertex). We will restrict our attention to finite connected ray-graphs.  We give some simple examples of models for $\cx$ in sections~\ref{sec:examples} and~\ref{sec:joins}.  \\

The metric on these graphs will be that of arc-length, and we will consider all edges as having length one.  We shall call the class of all such ray-graphs $\graphs$ and elements of that class $X$.  \\

We will sometimes refer to the \emph{containment hyperspace} $C_{A}(X) = \{ B \in C(X): A \subset B\}$.  This concept is especially useful to us when $A$ is a vertex of the graph.   If $A = \{p\}$ we may write $C_{p}(X)$ rather than $C_{\{p\}}(X)$.  \\

\section{Some Basic Examples of $(\cx, \haus)$}
\label{sec:examples}

In this section and the next we will discuss a few known models of the hyperspace $(\cx, \haus)$ for specific $X \in \graphs$.  For more detail on the hyperspaces in the first two subsections, see \cite{IllanesNadler}.

\subsection{$X \homeo [0,1]$}
\label{sec:interval}

\begin{figure}[h]
\begin{center}
\begin{tikzpicture}[scale=3]

\tikzstyle{shaded}=[fill=red!10!blue!20!gray!30!white]

\coordinate [label=below left:{$0$}] (p) at (-2, .5);
\coordinate [label=below right:{$1$}] (q) at (-1,.5);
\coordinate [label=below:{$X$}] (X) at (-1.5, .25);

\draw[gray] (p) -- (q);

\coordinate [label=below:{$(0,0)$}] (v1) at (0, 0);
\coordinate [label=above:{$(0,1)$}] (v2) at (0, 1);
\coordinate [label=above right:{$(1,1)$}] (v3) at (1,1);
\coordinate [label = right:{$C(X)$}] (CX) at (.5, .25);

\draw[gray] (v1) -- (v2);
\draw[gray] (v2) -- (v3);
\draw[gray] (v3) -- (v1);

\filldraw [shaded] (0,0)--(0,1)--(1,1)--(0,0);

\draw[thick] (-1.75, .5)--(-1.5, .5);
\coordinate [label = above:{$A$}] (A) at (-1.6, .5);

\coordinate [label=above:{$A$}] (pA) at (.25, .5);
\fill[black] (pA) circle (.4pt);

\end{tikzpicture}
\caption{$X = [0,1]$ and $\cx$, as well as an element $A \in \cx$}
\label{fig:triangle}
\end{center}
\end{figure}
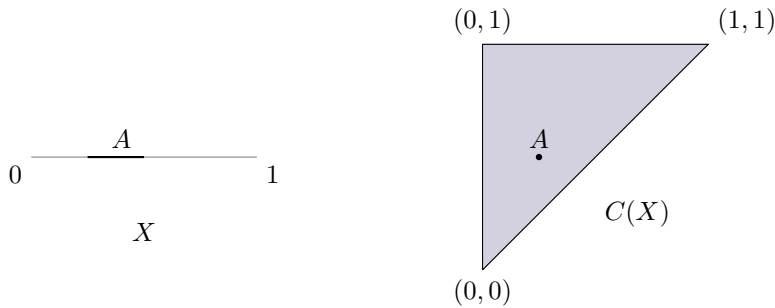

If $X$ is a segment homeomorphic to a closed interval then any element $A \in \cx$ is of the form $[a, b]$. Let us assume $X = [0,1]$, and then we have $0 \leq a \leq b \leq 1$.  There is a homeomorphism from the hyperspace $\cx$ to the solid triangle in $\R^{2}$ with vertices at $(0,0)$, $(0,1)$ and $(1,1)$ which takes an interval $[a, b]$ to the point $(a, b)$.  (Here we are abusing the notation to say that $[a, a] = \{a\}$.)  See Figure~\ref{fig:triangle}.  Notice that the left edge of the triangle corresponds to subsets of $X$ which contain $0$, i.e. the containment hyperspace $C_{\{0\}}(X)$.  The top edge corresponds to subsets which contain $1$, $C_{\{1\}}(X)$, and the hypotenuse corresponds to the single-element sets $F_{1}(X)$.  We will refer to this triangle as $T$.

\subsection{$X \homeo S^{1}$}
\label{sec:circle}

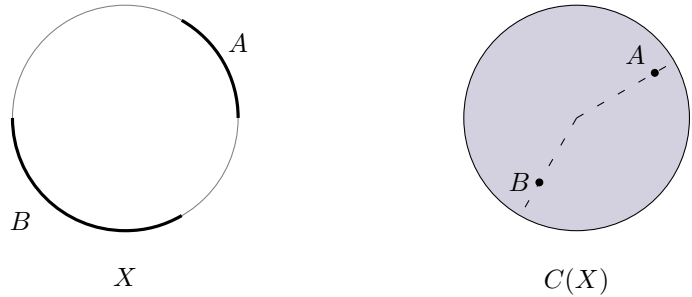
\begin{figure}[h]
\begin{center}
\begin{tikzpicture}[scale=1.5]

\tikzstyle{shaded}=[fill=red!10!blue!20!gray!30!white]

\draw[gray] (-4,0) circle (1cm);
\coordinate [label = below:{$X$}] (X) at (-4, -1.25);
\draw[very thick] (-3,0) arc (0:60:1cm);
\coordinate [label=above:{$A$}] (A) at (-3, .5);
\draw[very thick] (-5, 0) arc (180:300:1cm);
\coordinate [label=below left:{$B$}] (B) at (-4.75, -.75);

\filldraw [shaded] (0,0) circle (1cm);
\coordinate [label=below:{$C(X)$}] (CX) at (0, -1.25);
\coordinate [label=above left:{$A$}] (pA) at (canvas polar cs:angle=30,radius=.8cm);
\fill[black] (pA) circle (1pt);
\draw[loosely dashed] (0,0) -- (canvas polar cs:angle=30, radius=1cm);

\coordinate [label=left:{$B$}] (pB) at (canvas polar cs:angle=240, radius=.66cm);
\fill[black] (pB) circle (1pt);
\draw[loosely dashed] (0,0) -- (canvas polar cs:angle=240, radius=1cm);

\end{tikzpicture}
\caption{$X = S^{1}$ and $\cx$, as well as two elements $A, B \in \cx$}
\label{fig:disc}
\end{center}
\end{figure}

If $X$ is a simple closed curve, then elements of $\cx$ can be categorized by their midpoint and their length.  Let $X = S^{1}$.  We can make a homeomorphism from $\cx$ to the solid unit disc by mapping an arc with length $l$ and midpoint $p$ to the point which sits on the radial line through $p$, and whose distance from the origin is $1-\frac{l}{2\pi}$.  See Figure~\ref{fig:disc}.  Notice that the boundary of the disc corresponds to the single-element sets $F_{1}(X)$, and the center point of the disc corresponds to the full circle.  We will refer to this disc as $D$.\\

Although $[0,1] \not\homeo S^{1}$, their two hyperspaces $D$ and $T$ are homeomorphic.  It is known that for finite graphs this is the only such example \cite{Duda1}.  

\subsection{$X \homeo [0, \infty)$}
\label{sec:ray}

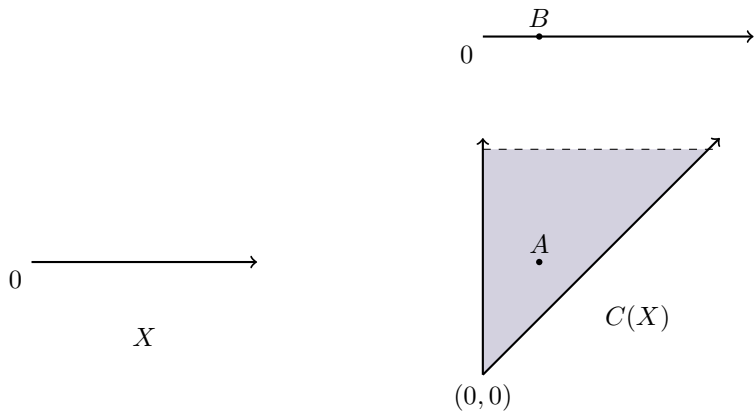
\begin{figure}[h]
\begin{center}

\begin{tikzpicture}[scale=3]
\tikzstyle{shaded}=[fill=red!10!blue!20!gray!30!white]

\coordinate [label = below left:{$0$}] (p) at (-2, .5);
\coordinate (q) at (-1, .5);
\coordinate [label=below:{$X$}] (X) at (-1.5, .25);
\draw[->, thick] (p) -- (q);

%

\coordinate [label=below:{$(0,0)$}] (v1) at (0, 0);
\coordinate (v2) at (0, 1.05);
\coordinate (v3) at (1.05,1.05);
\coordinate [label = right:{$C(X)$}] (CX) at (.5, .25);
\fill [shaded] (0,0)--(0,1)--(1,1);
\draw[->, thick] (v1) -- (v2);
\draw[->, thick] (v1) -- (v3);
\draw[dashed] (0,1) -- (1.02,1); 

\coordinate [label= below left:{$0$}] (p1) at (0, 1.5);
\coordinate (p2) at (1.2, 1.5);
\draw[->, thick] (p1)--(p2);

\coordinate [label = above:{$A$}] (pA) at (.25, .5);
\fill[black] (pA) circle (.4pt);
\coordinate [label = above:{$B$}] (pB) at (.25, 1.5);
\fill[black] (pB) circle (.4pt);

\end{tikzpicture}
\caption{$X = [0, \infty)$ and $\cx$, as well as the elements $A = [.25, .5]$ and $B = [.25, \infty)$, both in $\cx$}
\label{fig:rayspace}
\end{center}
\end{figure}

If $X$ is a ray homeomorphic to $[0, \infty)$, then elements of $\cx$ are either bounded intervals of the form $[a, b]$ or unbounded intervals of the form $[a, \infty)$.  We can make a homeomorphism from $\cx$ to the space $T^{\infty} \sqcup [0, \infty)$, where $T^{\infty} = \{(a, b) \in \R^{2}: 0 \leq a \leq b \}$ is an ``infinite triangle."  This is done by mapping $[a, b]$ to $(a, b) \in T^{\infty}$ and $[a, \infty)$ to $a \in [0, \infty)$.  See Figure~\ref{fig:rayspace}.  \\

Notice that although $T^{\infty}$ is itself unbounded, elements of $T^{\infty}$ correspond to bounded subsets of $X$, and in particular, the left edge of $T^{\infty}$ corresponds to bounded elements which contain $0$, and the hypotenuse corresponds to single-element sets.  For a fixed horizontal value $a$, increasing the vertical value $b$ corresponds to longer bounded intervals.  Since the second component, $[0, \infty)$, corresponds to unbounded intervals, it can be loosely thought of as the ``top" of the infinite triangle. In this example, unlike before, the containment hyperspace $C_{\{0\}}(X)$ has two components: the left edge of the triangle \emph{and} the leftmost point $0 \in [0, \infty)$.  Clearly this $C(X)$ is not connected and not compact.\\


With these three examples we can form several more examples by understanding what happens to the hyperspace when you attach two graphs together in a specific way.

\section{The wedge product of graphs}
\label{sec:joins}

We begin by describing the hyperspace of the compound ray-graph $X = X_{1} \vee_{p} X_{2}$, where $p$ is a vertex of both, when models for the two hyperspaces $C(X_{1})$ and $C(X_{2})$ are already given. \\

It is clear that the hyperspace $\cx$ will contain all the elements which are in $C(X_{1})$ and $C(X_{2})$.  It will also contain elements which correspond to subsets of $X$ that contain the joining point $p$ and part of $X_{1}$ and $X_{2}$.  In fact, to any subset $A \subset X_{1}$ which contains $p$, we can union a subset of $X_{2}$ which contains $p$, and arrive at an element of $\cx$.  This shows that $\cx$ will contain a cross product of $\cpx{1}$ and $\cpx{2}$.

\begin{theorem}
\label{thm:join}
For $X_{1}, X_{2} \in \graphs$ and the wedge produce $X = X_{1} \vee_{p} X_{2}$ (where $p$ is a vertex of $X_{1}$ and $X_{2}$), then 
$$C(X) \homeo \frac{C(X_{1}) \sqcup \cpx{1} \times \cpx{2} \sqcup C(X_{2})}{(\cpx{1} \sim \cpx{1} \times \{p\} \mbox{ and } \{p\} \times \cpx{2} \sim \cpx{2})}$$ 
\end{theorem}

This theorem, which seems to be well-known in the folklore (and certainly applies to a larger class of spaces than graphs), gives us the following nice algorithm for drawing $\cx$:

\begin{enumerate}
\item Draw $\cpx{1} \times \cpx{2}$.
\item Attach the rest of $C(X_{1})$ to the figure by identifying its subset $\cpx{1}$ with the slice $\cpx{1} \times \{p\}$ in the cross product.  
\item Attach the rest of $C(X_{2})$ to the figure by identifying its subset $\cpx{2}$ with the slice $\{p\} \times \cpx{2}$ in the cross product.
\end{enumerate}

We shall use this algorithm in the following examples.

\subsection{The noose}
\label{sec:noose}

We begin with a simple example, the noose.  The space is the wedge-product of a circle and an interval.  To follow the steps outlined in the algorithm, let $X_{1} \homeo S^{1}$ and $X_{2} \homeo [0,1]$, and let $p$ be the point 0 on both.  Begin by noting that $\cpx{1}$ is homeomorphic to a subdisc inside the disc $D$, which includes the point $p$.  (In fact, the shape is a cardiod, but for ease of representation we will draw it as a disc.)  $\cpx{2}$ is the left edge of the triangle $T$.  \\

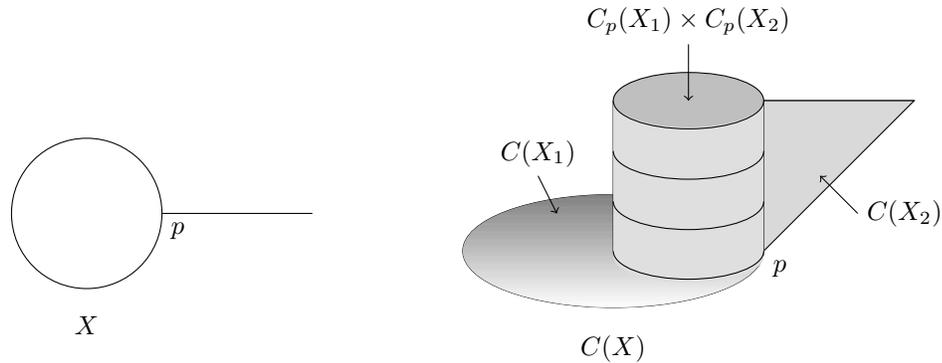
\begin{figure}[h]
\begin{center}
\begin{tikzpicture}[scale=1]
\tikzstyle{shaded}=[fill=red!10!blue!20!gray!30!white]

\coordinate [label=below right:{$p$}] (p1) at (-6,.5);
\draw (-7,.5) circle (1cm);
\draw (p1) -- (-4, .5);
\coordinate [label=below:{$X$}] (X) at (-7, -.75);

\coordinate (c1) at (0,0);
\coordinate (c2) at (1, 0);
\coordinate [label=below right:{$p$}] (p) at (2, 0);
\coordinate (q1) at (0,2);
\coordinate (q2) at (2,2);
\coordinate(q3) at (4, 2);

\draw (0,0) ellipse (2cm and .75cm);
\shade[bottom color = white, top color = gray] (0,0) ellipse (2cm and .75cm);
\coordinate [label = above:{$C(X_{1})$}] (CX1) at (-1,1);
\coordinate (x1) at (-.75, .5);
\draw[->] (CX1) -- (x1);

\draw (c1)--(q1);
\draw (q2) -- (q3);
\fill[gray!30] (p) -- (q2) -- (q3) -- (p);
\draw (p) -- (q2) -- (q3) -- (p);
\coordinate [label= right:{$C(X_{2})$}] (CX2) at (3.25, .5);
\coordinate (x2) at (2.75, 1);
\draw[->] (CX2) -- (x2);

\node [cylinder, rotate = 90, minimum height = 2.8cm, minimum width = 2cm, aspect=3.1, cylinder uses custom fill, cylinder end fill=gray!50, cylinder body fill=gray!25] at (1,.7) {};
\draw (0,0) arc (180:360:1cm and .375cm);
\draw (1,2) ellipse (1cm and .375cm);
\draw (0,.66) arc (180:360:1cm and .375cm);
\draw (0,1.33) arc (180:360:1cm and .375cm);
\coordinate [label = above:{$\cpx{1} \times \cpx{2}$}] (CP) at (1,2.75);
\coordinate (x3) at (1,2);
\draw[->] (CP) --(x3);

\coordinate [label=below:{$C(X)$}] (CX) at (0, -1);

\end{tikzpicture}
\caption{$X$ a noose, and $\cx$}
\label{fig:noose}
\end{center}
\end{figure}

When we cross $\cpx{1}$ with $\cpx{2}$ we get a solid cylinder.  Then we attach the rest of the disc $D = C(X_{1})$ along the slice $\cpx{1} \times \{p\}$, which is the bottom of the cylinder.  Finally, we attach the rest of the triangle $T = C(X_{2})$ along the slice $\{p\} \times \cpx{2}$, which is a vertical line on the boundary of the cylinder.  Note that the resulting space has both two- and three-dimensional sections.  See Figure~\ref{fig:noose}.

\subsection{$n$-od}
\label{sec:nod}

An $n$-od is a point $p$ with $n$ intervals attached such that they intersect only at $p$.  We begin with the $2$-od.  Of course the $2$-od is homeomorphically the same as an interval, and so we should arrive at $C(X) \homeo T$.  But for the purposes of the construction, we will go through it nonetheless.  In this case $X_{1} = X_{2} \homeo [0,1]$ and $\cpx{1} = \cpx{2}$ is the left edge of the triangle $T$.  \\

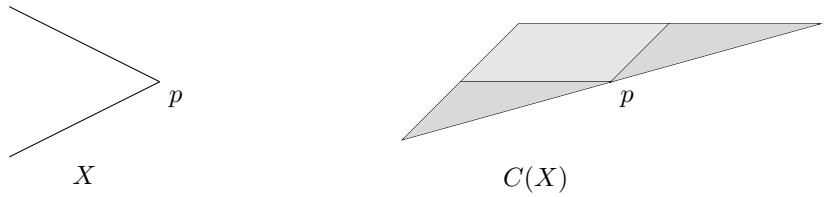
\begin{figure}[h]
\begin{center}

\begin{tikzpicture}[scale=2]
\tikzstyle{shaded}=[fill=red!10!blue!20!gray!30!white]

\coordinate [label = below right:{$p$}] (p) at (-2, 0, 0);
\coordinate (q1) at (-3, .5, 0);
\coordinate (q3) at (-3, -.5, 0);

\draw (p)--(q1);
\draw (p)--(q3);

\coordinate [label=below:{$X$}] (X) at (-2.5, -.5, 0);

\draw (0,0,0)--(0,0,-1)--(1,0,-1)--(1,0,0)--(0,0,0);
\fill[gray!20] (0,0,0)--(0,0,-1)--(1,0,-1)--(1,0,0)--(0,0,0);

\coordinate [label=below right:{$p$}] (p1) at (1,0,0);

\draw (0,0,0)--(0,0,1)--(1,0,0)--(0,0,0);
\fill[gray!30] (0,0,0)--(0,0,1)--(1,0,0)--(0,0,0);

\draw (1,0,0)--(1,0,-1)--(2,0,-1)--(1,0,0);
\fill[gray!30] (1,0,0)--(1,0,-1)--(2,0,-1)--(1,0,0);

\coordinate [label = below:{$C(X)$}] (CX) at (.5, -.5, 0);

\end{tikzpicture}

\caption{$X$ a 2-od and $C(X)$}
\label{fig:2od}
\end{center}
\end{figure}

When we cross $\cpx{1}$ with $\cpx{2}$ we get a square; we then attach the triangles $C(X_{i})$ along the appropriate edges.  The result is, of course, another triangle.  See Figure~\ref{fig:2od}. \\

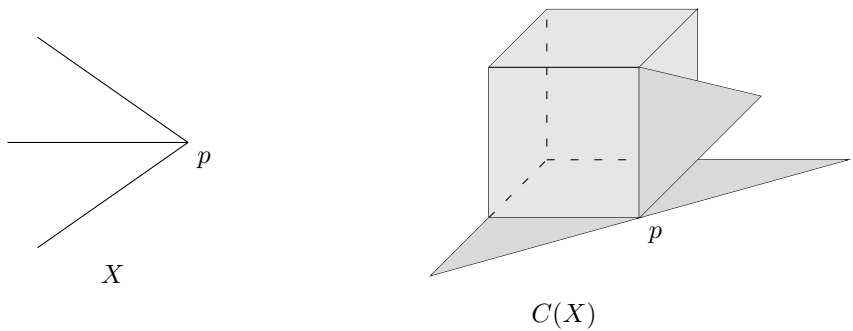
\begin{figure}[h]
\begin{center}

\begin{tikzpicture}[scale=2]
\tikzstyle{shaded}=[fill=red!10!blue!20!gray!30!white]

\coordinate [label = below right:{$p$}] (p) at (-2, .5, 0);
\coordinate (q1) at (-3, 1.2, 0);
\coordinate (q2) at (-3.2, .5, 0);
\coordinate (q3) at (-3, -.2, 0);

\draw (p)--(q1);
\draw (p)--(q2);
\draw (p)--(q3);

\coordinate [label=below:{$X$}] (X) at (-2.5, -.25, 0);

\draw (0,0,0)--(1,0,0)--(1,1,0)--(0,1,0)--(0,0,0);
\fill[gray!20] (0,0,0)--(1,0,0)--(1,1,0)--(0,1,0)--(0,0,0);

\draw (0,1,0)--(0,1,-1)--(1,1,-1)--(1,1,0);
\fill[gray!20] (0,1,0)--(0,1,-1)--(1,1,-1)--(1,1,0);

\draw (1,0,0)--(1,0,-1)--(1,1,-1)--(1,1,0)--(1,0,0);
\fill[gray!20] (1,0,0)--(1,0,-1)--(1,1,-1)--(1,1,0)--(1,0,0);

\draw[loosely dashed] (0,0,0)--(0,0,-1)--(0,1,-1);
\draw[loosely dashed] (0,0,-1)--(1,0,-1);
\draw (0, 1, 0)--(1, 1, 0);

\coordinate [label=below right:{$p$}] (p1) at (1,0,0);

\draw (0,0,0)--(0,0,1)--(1,0,0)--(0,0,0);
\fill[gray!30] (0,0,0)--(0,0,1)--(1,0,0)--(0,0,0);

\draw (1,0,0)--(1,0,-1)--(2,0,-1)--(1,0,0);
\fill[gray!30] (1,0,0)--(1,0,-1)--(2,0,-1)--(1,0,0);

\draw (1,0,0)--(1,1,0)--(2,1, .5)--(1,0,0);
\fill[gray!30] (1,0,0)--(1,1,0)--(2,1,.5)--(1,0,0);

\coordinate [label = below:{$C(X)$}] (CX) at (.5, -.5, 0);

\end{tikzpicture}
\caption{$X$ a triod and $\cx$, a cube with three fins.}
\label{fig:triod}
\end{center}
\end{figure}

Now we will use the $2$-od to construct the hyperspace of the triod, by letting $X_{1}$ be the 2-od and $X_{2}$ be the interval.  $\cpx{1}$ is now the square from before, and $\cpx{2}$ is again the left edge of $T$.  The cross product therefore results in a cube.  The two fins from the hyperspace of the 2-od are attached along the bottom of the cube, and another fin (the rest of $C(X_{2})$) is attached along the front right edge.  See Figure~\ref{fig:triod}. \\

It is easy to see that continuing in this manner will result in the hyperspace of the $n$-od consisting of the $n$-cube $I_{n}$, with $n$ fins attached along its edges, each with a corner at $p$.\\

\subsection{The infinite noose}
\label{sec:infinitenoose}

Let $X$ be the infinite noose, made up of $X_{1} \homeo S^{1}$ and $X_{2} \homeo [0, \infty)$, joined at the point $p$.  $\cpx{1}$, as we have already noted in section~\ref{sec:noose}, is a cardiod inside the unit disc, which we will draw as a subdisc.  Recall from section~\ref{sec:ray} that $\cpx{2}$ has two components: the left edge of the infinite triangle $T^{\infty}$, and the left-most point of the ray.  The point of $C(X_{2})$ which corresponds to the single-point set $\{p\} \subset X_{2}$ is at the bottom of the triangle.  \\

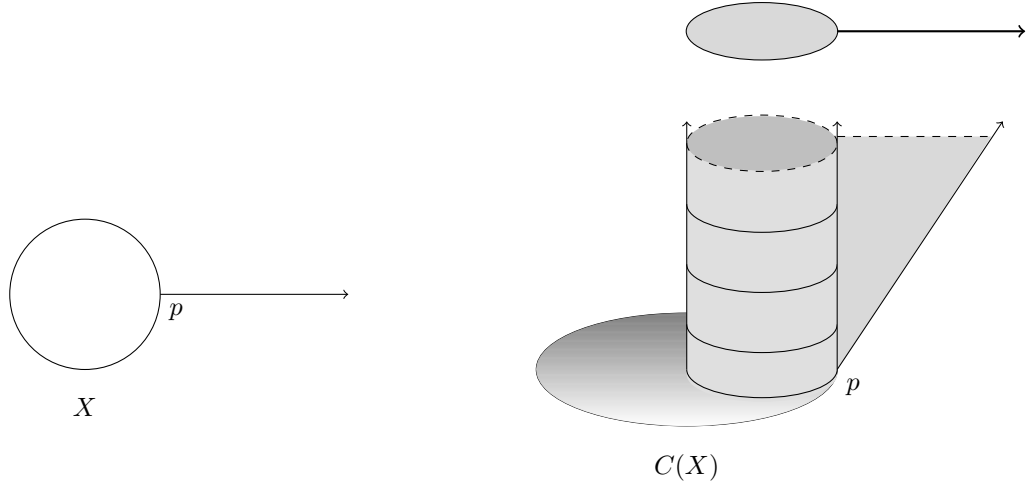
\begin{figure}[h]
\begin{center}

\begin{tikzpicture}[scale=1]
\tikzstyle{shaded}=[fill=red!10!blue!20!gray!30!white]

\coordinate [label=below right:{$p$}] (p1) at (-7,1);
\draw (-8,1) circle (1cm);
\draw[->] (p1) -- (-4.5, 1);
\coordinate [label=below:{$X$}] (X) at (-8, -.25);

\coordinate (c1) at (0,0);
\coordinate (c2) at (1, 0);
\coordinate [label=below right:{$p$}] (p) at (2, 0);
\coordinate (q1) at (0,3.1);
\coordinate (q2) at (2,3.1);
\coordinate(q3) at (4, 3.1);

\draw (0,0) ellipse (2cm and .75cm);
\shade[bottom color = white, top color = gray] (0,0) ellipse (2cm and .75cm);

\node [cylinder, rotate = 90, minimum height = 3.8cm, minimum width = 2cm, aspect=3.1, 
 cylinder uses custom fill, cylinder end fill=gray!50, cylinder body fill=gray!25] at (1,1.2) {};
\draw (0,0) arc (180:360:1cm and .375cm);
\draw (0,.6) arc (180:360:1cm and .375cm);
\draw (0,1.4) arc (180:360:1cm and .375cm);
\draw (0,2.2) arc (180:360:1cm and .375cm);
\draw[dashed] (0,3.01) arc (180:360:1cm and .375cm);
\draw[dashed] (2,3.01) arc (0:180:1cm and .375cm);

\fill[gray!30] (p) -- (q2) -- (4.05, 3.1) -- (p);
\draw[->] (c1)--(0,3.3);
\draw[->] (p) -- (2,3.3);
\draw[->] (p) -- (4.2, 3.3);
\draw[dashed] (q2) -- (q3);

\draw[thick] (1,4.5) ellipse (1cm and .375cm);
\fill[gray!30] (1,4.5) ellipse (1cm and .375cm);
\draw[->, thick] (2, 4.5) -- (4.5, 4.5);

\coordinate [label=below:{$C(X)$}] (CX) at (0, -1);

\end{tikzpicture}
\caption{$X$ is the infinite noose. $\cx$ has two components.}
\label{fig:infinnoose}
\end{center}
\end{figure}

Crossing $\cpx{1}$ with $\cpx{2}$, we get an infinite cylinder and a disc.  We attach $C(X_{1})$ to the slice $\cpx{1} \times \{p\}$, along the bottom of the cylinder.  We attach $C(X_{2})$ along $\{p\} \times \cpx{2}$, producing an infinite fin off the side of the cylinder and a ray off the side of the disc.  See Figure~\ref{fig:infinnoose}.

\subsection{The real line}
\label{sec:realline}

Let $X_{1} = X_{2} \homeo [0, \infty)$, both just a single vertex with a ray attached.  Then we can think of $X \homeo \R$ as the result of attaching these two subgraphs along their vertex.  Both subgraphs have hyperspaces which consist of $T^{\infty} \sqcup [0, \infty)$, and the containment hyperspace for the vertex is the union of the left edge of the triangle and the leftmost point of the ray. \\

Following the algorithm, $\cpx{1} \times \cpx{2}$ gives us four components: an infinite square, two rays, and a point.   When we attach the rest of $C(X_{1})$ along the correct slice, it attaches the rest of the triangle along one side of the infinite square, and the rest of a ray along one of the rays.  Similarly when we attach the rest of $C(X_{2})$ it attaches the rest of an infinite triangle along the other side of the infinite square, and another ray along the second ray.  The end result is four components: a half-plane, two real lines, and a point. \\

To see this more algebraically, we will briefly construct a homeomorphism to show:
$$(\cx, \haus) \homeo \{0\} \sqcup \R \sqcup \R \sqcup (\R \times [0, \infty))$$

We will write different sections of the disjoint union with a subscript to distinguish them, e.g. $\R_{1}$ and $\R_{2}$.  Because elements of $\cx$ are connected closed subsets of $X$, they are closed intervals, and are therefore one of four types: $A = X = \R$, $A = [a, \infty)$ (``unbounded to the right"),  $A = (-\infty, a]$ (``unbounded to the left") or $A = [a, b]$ (``bounded"). \\

Define a map $\phi$ from $\cx \rightarrow \{0\} \sqcup \R_{1} \sqcup \R_{2} \sqcup (\R \times [0, \infty))$ in the following way:
 \[ 
 \phi(A) = \left\{ \begin{array}{ll}
 	\{0\} & \mbox{ if $A = \R$} \\
	a \in \R_{1} & \mbox{ if $A = [a, \infty)$ } \\
	a \in \R_{2} & \mbox{ if $A = (-\infty, a]$ } \\
	(a, b-a) \in \R \times [0, \infty) & \mbox{ if $A = [a, b]$ } 
	\end{array}
	\right. \]

This is clearly a homeomorphism.  Notice that the four components of the hyperspace correspond to the four different ways in which it is possible for a subset to be unbounded.  In section~\ref{sec:hausdorff} we will show this is not a coincidence.

\section{The connected components of $(\cnx, \haus)$}
\label{sec:hausdorff}

The last two examples of section~\ref{sec:joins} show that under the Hausdorff topology, there is a relationship between the number of rays in a given graph, and the number of connected components of its hyperspace.  That relationship is what we explore in this section. \\

We begin by developing some extra terminology to deal with ray-graphs in particular.  Let $\mathcal{R} = \{R_{1}, \dots R_{k}\}$ denote the set of rays in a given ray-graph.  If \#$\mathcal{R} = k$, we will call $X$ an \emph{$k$-legged} graph.  We will denote by $X_{G} = X - \cup_{i=1}^{k} R_{i}$.  If $A \subset X$, and $A \cap R_{i}$ is an unbounded interval, we say that $A$ is \emph{unbounded in direction $i$}.  In this way we can talk about the \emph{unbounded direction set of $A$}, which is the set of indices between 1 and $k$ for which $A$ is unbounded in direction $i$.  Clearly there are $2^{k}$ possible unbounded direction sets, in one-to-one correspondence with the power set of $\{1, 2, \dots, k\}$.\\

Let $\mathcal{P}_{k}$ be the power set of $\{1, 2, \dots, k\}$.  Define a function $\phi: \cnx \rightarrow \mathcal{P}_{k}$ by $\phi(A) = \Delta$, where $\Delta \in \mathcal{P}_{k}$ is the unbounded direction set of $A$. Recall that we denote the Hausdorff distance between two elements $A$ and $B$ by $d_{H}(A, B)$. \\

\begin{lemma}
Let $A, B \in \cnx$ under the Hausdorff topology.
\begin{enumerate}
\item If $d_{H}(A, B) < \infty$, then $A$ and $B$ have the same unbounded direction set.
\item If $A$ and $B$ have distinct unbounded direction sets, e.g. there exists a ray $R_{i} \in X$ such that $A$ is unbounded in direction $i$ but $B$ is not, then there does not exist any path through $\cnx$ from $A$ to $B$.
\end{enumerate}
\label{lem:1}
\end{lemma}

\textbf{Proof:}
If $A$ is unbounded in direction $i$ and $B$ is not, then clearly $d_{H}(A, B) = \infty$.  Since any path is a continuous image of a compact set, it must have a compact image which contains $A$ and $B$.  If $d_{H}(A, B) = \infty$, this is impossible. \hfill $\Box$ \\
 
\begin{lemma}
If $X \in \graphs$ is a $k$-legged graph, then for all $n$, the hyperspace $(\cnx, \haus)$ has at least $2^{k}$ connected components.
\label{lem:part1}
\end{lemma}

\textbf{Proof:} Let $X \in \graphs$ be a graph with $k$ distinct rays, labelled $R_{1}, \dots, R_{k}$.  Consider the power set $\mathcal{P}_{k}$ of $\{1, \dots, k\}$.  Each element of the power set corresponds to an unbounded direction set.   We will use the map $\phi: \cnx \rightarrow \mathcal{P}_{k}$ from above, given by $\phi(A) = \Delta$ if $A$ is unbounded in the direction set $\Delta$.  We will show that $\phi$ is continuous, and therefore $\cnx$ has at least $2^{k}$ connected components. \\

Because $(\cnx, \haus)$ is first countable, it is enough to show convergent sequences are mapped to convergent sequences.  Let $A_{m} \rightarrow A$ be a convergent sequence of elements of $\cnx$, meaning that $d_{H}(A_{m}, A) \rightarrow 0$ as $m \rightarrow \infty$.  If $A$ is unbounded in direction $R_{i}$, and $A_{m}$ is not (or vice versa) we know $d_{H}(A, A_{m}) = \infty$, so for all $m$ greater than some $m^{*}$ we must have $A_{m}$ unbounded in the same set of directions as $A$.  Therefore $\phi(A_{m}) = \phi(A)$ for all $m > m^{*}$ and $\phi$ is continuous.
\hfill $\Box$ 

\begin{theorem}
If $X \in \graphs$ is a $k$-legged graph, then for all $n$, the hyperspace $(\cnx, \haus)$ has \textbf{exactly} $2^{k}$ path-connected components.
\label{thm:hausdorffcomponents}
\end{theorem}

The previous lemma showed that $\cnx$ has at least $2^{k}$ connected components.  We will now show that it has no more than that, by showing that for all $\Delta \in \mathcal{P}_{k}$, $\{A \in \cnx: \phi(A) = \Delta\}$ is a path-connected set.   This will be done by taking any element in a given component and constructing a path from it to a designated ``default" element of that component.  \\

A note on notation:  a subinterval of $A$ inside the edge $E_{i}$ will be denoted $[a,b]_{i}^{e}$.  A subinterval of $A$ inside the ray $R_{i}$ will be denoted $[a, b]_{i}^{r}$. For rays, the vertex is 0.  For edges which only have one ramification point $X$, the endpoint 0 is the ramification point.   \\

\textbf{Proof:} 
Fix $n$.  We begin by choosing for each $\Delta \in \mathcal{P}_{k}$ a particular element $A_{\Delta}$ of $\{A \in \cnx: \phi(A) = \Delta\}$.  The element $A_{\Delta}$ will consist of the complete finite-graph $X_{G}$, and all the rays which are in the unbounded direction set $\Delta$, but no part of the other rays.  It will have one connected component.  Precisely, 

$$A_{\Delta} = X_{G} \cup \bigcup_{i \in \Delta} R_{i}$$
  
Given an element $A \in \cnx$ with $\phi(A) = \Delta$, we will construct a path from $A$ to $A_{\Delta}$.  There are three steps.  First, any sections of $A$ contained completely in rays $R_{i}$ where $i \in \Delta$ we will grow until they touch $X_{G}$ at the vertex.  Secondly, any sections of $A$ contained in rays $R_{j}$ where $j \not\in \Delta$, we will shrink down until they are gone.  Finally, we grow the remaining subset out so that it includes all of $X_{G}$.  The first two steps of this process will either keep constant or decrease the number of components of $A$; the last step will produce an element with one component.  So the path will stay in $\cnx$ at all times.  \\

By definition, if $i \in \Delta$ then $A \cap R_{i} \neq \emptyset$.  Consider those $i \in \Delta$ for which $A \cap R_{i} \neq R_{i}$.  For each such $i$, that intersection will be a finite number of intervals, one of which is unbounded.  Call the unbounded one $[a_{i}, \infty)_{i}^{r}$, and call the vertex where that ray is attached $v_{i}$.  We will grow this interval out so that it encompasses all of $R_{i}$.   We define the first step in the path as follows.  $f_{0}:[0,1] \rightarrow \cnx$ is given by

$$f_{0}(t) = A \cup \bigcup_{i \in \Delta} [tv_{i} + (1-t)a_{i}, \infty)_{i}^{r}$$

If $A$ had several intervals contained in that ray, this process will consume them, reducing the number of components of $A$.  If $A$ had empty unbounded direction set, this will do nothing to $A$.  Let $A_{1} = f_{0}(1)$. \\ 

Now consider those $j \not\in \Delta$ with $A_{1} \cap R_{j} \neq \emptyset$.  That intersection will consist of a finite number of bounded intervals $[a_{j}^{i}, b_{j}^{i}]_{j}^{r}$ (where $i = 1, \dots, l_{j}$ for some $l_{j} \leq k$).  We wish to shrink and slide each of those intersections down to the vertex $v_{j}$.  To do that we define a path $f_{1}:[0,1] \rightarrow \cnx$ by

$$f_{1}(t) = A_{1}  \cup \bigcup_{j \not\in \Delta} \bigcup_{i = 1}^{l_{j}} [tv_{j} + (1-t)a_{j}^{i}, tv_{j} + (t-1)b_{j}^{i}]_{j}^{r}$$

This is a path from $A_{1}$ to the set which agrees with $A_{1}$ in $X_{G}$, contains all of the rays in the unbounded direction set, but does not contain any section of any rays which are not in the unbounded direction set (apart from possibly the vertices).  Call this second intermediate set $A_{2} = f_{1}(1)$. \\

The final step will grow the subset $A_{2}$ out until it includes all of $X_{G}$.  Fix an element $a \in A_{2} \cap X_{G}$.  Because $X_{G}$ is a graph, it is path connected, so there exists a path $\gamma:[0, 1] \rightarrow X_{G}$ which starts at $a$ and whose image contains all of $X_{G}$.  Define

$$f_{2}(t) = A_{2} \cup \bigcup_{x \in [0, t]} \gamma(x)$$

Clearly $f_{2}(0) = A_{2}$ and $f_{2}(1) = A_{\Delta}$.  The continuity of $\gamma$ makes $f_{2}$ continuous, and because components may merge together, but never split apart, the construction ensures $f_{2}(t) \in \cnx$ at all times.  Following $f_{0}$ with $f_{1}$ and $f_{2}$, we have a path from $A$ to $A_{\Delta}$.  Hence the set $\{A \in \cnx: \phi(A) = \Delta\}$ is path-connected.  This completes the proof.
\hfill $\Box$ \\

\section{Connectedness of $(\cnx, \vietoris)$}
\label{sec:vietoris}

In this section we will explore some of the distinctions between a hyperspace of a finite ray-graph under the Hausdorff topology with that same hyperspace under the Vietoris topology.  We begin by recalling the definition. 

\subsection{The Vietoris topology}

Let $X$ be the base space, and let $U_{1}, \dots, U_{n}$ be a finite number of open subsets of $X$.  For any hyperspace $\Hx$ over $X$, we define the open set $U^{*} = <U_{1}, \dots, U_{n}>$ in the following way: 

$A \in U^{*}$ iff
\begin{enumerate}
\item $A \subset \bigcup_{i=1}^{n} U_{i}$
\item $A \cap U_{i} \neq \emptyset$ for all $i = 1, \dots, n$
\end{enumerate}

Such open sets make a basis for the Vietoris topology on $\Hx$.  It is sometimes more useful to treat the topology as the supremum of the upper and lower Vietoris topologies.  The upper Vietoris topology is generated by sets of the form 
$$U^{+} = \{ A \in \Hx : A \subset U \}$$

where $U$ is open in $X$.  The lower Vietoris topology is generated by sets of the form 
$$V^{-} = \{A \in \Hx: A \cap V \neq \emptyset\}$$

where $V$ is open in $X$. Subbase elements of the Vietoris topology are then of the form $U^{+}
$ and $V^{-}_{1} \cap V^{-}_{2} \cap \dots \cap V_{n}^{-}$.  

\subsection{Path-connectedness}
\label{sec:connected}

In~\cite{Esty} we proved that $CL(M)$ was contractible for any Borel compact space $M$ having the property that the closure of open balls is closed balls.  Since ray-graphs satisfy those conditions, we know that $\clx$ is contractible:

\begin{theorem}
$(CL(X), \vietoris)$ is contractible.
\end{theorem}

We now prove a companion theorem to Theorem~\ref{thm:hausdorffcomponents}.

\begin{theorem} 
\label{thm:cx}
$(\cnx, \vietoris)$ is path-connected. 
\end{theorem}

The proof will be by construction and is similar in flavor to the proof of Theorem~\ref{thm:hausdorffcomponents}.  In fact, the first part is identical: take an element $A \in \cx$ and construct a path from it to the element $A_{\Delta}$, $Delta = \phi(A)$. As it was in the Hausdorff topology, this construction is continuous in the Vietoris topology.  The distinction comes when we then form a path from $A_{\Delta}$ to the element $X$.  \\

\textbf{Proof:}  Recall that if $\phi(A) = \Delta \in \mathcal{P}_{k}$, we define the element $A_{\Delta}$ as 
$$ A_{\Delta} = X_{G} \cup \bigcup_{i \in \Delta} R_{i}$$

Start with the same path $f = f_{2} \circ f_{1} \circ f_{0}$ from $A$ to $A_{\Delta}$ as given in the proof of Theorem~\ref{thm:hausdorffcomponents}.  As before, the path remains in $\cnx$ at all times, and results in $A_{\Delta} \in \cx$.  The proof that this is continuous under the Vietoris topology is similar to the proof that the second path is continuous.  We will prove the latter.  \\

For the second path we will connect $A_{\Delta}$ to $X$.  Let $f(t) = \frac{t}{1-t}$, and define a path $\gamma: [0, 1] \rightarrow C(X)$ by:

\[ \gamma(t) = \left\{ \begin{array}{ll}
 A_{\Delta} \cup \bigcup_{i \not\in \Delta} [0, f(t)]_{i}^{r} & t \in [0, 1) \\
 X & t = 1 
 \end{array} \right. \]

Obviously $\gamma(0) = A_{\Delta}$ and by construction, $\gamma(t) \in \cx$ for all $t$.   To show $\gamma$ is continuous, it is enough to show that it is continuous with respect to the upper and lower Vietoris topologies.  Note that $s < t$ implies $\gamma(s) \subset \gamma(t)$.  This means for the upper topology, which is concerned with containment, we only have to worry about $t$ increasing, and for the lower topology, which is concerned with intersection, we only have to worry about $t$ decreasing.\\

We begin by checking continuity with respect to the upper Vietoris topology.  Fix $t_{0} \in [0,1]$ and suppose $\gamma(t_{0}) \in U^{+}$.  If $t_{0} = 1$ then $\gamma(t_{0}) = X$ and as $X \in U^{+}$ we have $U = X$, which clearly implies $\gamma(t) \in U^{+}$ for all $t$.  So we assume that $t_{0} \in [0, 1)$ and that $U \neq X$. \\

Since $U \neq X$ and $\gamma(t_{0}) \subset U$, there exists $\epsilon = d(\gamma(t_{0}), U^{c}) > 0$.  Continuity of $f(t)$ implies that there exists some $\delta$, $0 < \delta \leq \epsilon$ such that if $|t-t_{0}| < \delta$ then $|f(t) - f(t_{0})| < \epsilon$.  Then the total growth from $\gamma(t_{0})$ to $\gamma(t)$ is small, i.e. $d(\gamma(t), \gamma(t_{0})) < \epsilon$ and hence $\gamma(t) \in U^{+}$. So $\gamma$ is continuous with respect to the upper Vietoris topology. \\

Now we check continuity with respect to the lower Vietoris topology.  Let $\gamma(t_{0}) \in V_{1}^{-} \cap V_{2}^{-} \cap \dots \cap V_{n}^{-}$, meaning that $\gamma(t_{0}) \cap V_{i} \neq \emptyset$ for each $i$.  Pick a point $x_{i} \in \gamma(t_{0}) \cap V_{i}$.  If $x_{i} \in A_{\Delta}$ then $x_{i} \in \gamma(t)$ for all $t$.  If not, then $x_{i} \in [0, f(t_{0})]_{\ell}^{r}$ for some $\ell$.  Possibly $x_{i}$ is in the interior of that interval, or possibly it is the endpoint $f(t_{0})$.  If $x_{i}$ is in the interior, i.e. $x_{i} < f(t_{0})$, then pick $\delta_{i} > 0$ such that $|t-t_{0}| < \delta_{i}$ implies $|f(t_{0}) - f(t)| < f(t_{0}) - x_{i}$, and then $x_{i} \in f(t)$ also.  If $x_{i} = f(t_{0})$, then let $d_{i} = d(x_{i}, V_{i}^{c}) > 0$ and choose $\delta_{i}$ such that $|t-t_{0}| < \delta_{i}$ implies $|f(t) - f(t_{0})| < d_{i}$.  Then $[0, f(t)]_{\ell}^{r} \cap V_{i} \neq \emptyset$.  All together, let $\delta = \min\{ \delta_{i}: i = 1, \dots, n\}$.  Then for each $i$, $\gamma(t) \cap V_{i} \neq \emptyset$.  So $\gamma$ is continuous with respect to the lower Vietoris topology. \\

Combining with the path $f$, we have constructed a path from any $A \in \cnx$ to $X$.   \hfill $\Box$ \\

\bibliographystyle{plain}

\begin{thebibliography}{10}

\smallskip


\bibitem{Acosta1} G. Acosta, \textit{Continua with unique hyperspace},
Continuum theory; proceedings of the special session in honor of Professor Sam B. Nadler, Jr.'s 60 birthday. \textbf{230} (2002), 33-49.

\smallskip

\bibitem{Acosta2} G. Acosta, \textit{Continua with almost unique hyperspace},
Top. App.  \textbf{117} (2002), 175-189.

\smallskip


\bibitem{Duda1} R. Duda, \textit{On the hyperspace of subcontinua of a finite graph I},
Fundamenta Mathematicae. \textbf{62} (1968), 265--286.

\smallskip

\bibitem{Duda2} R. Duda, \textit{On the hyperspace of subcontinua of a finite graph II},
Fundamenta Mathematicae. \textbf{63} (1968), 225--255.

\smallskip

\bibitem{EberhartNadler} C. Eberhart,  S. Nadler, \textit{Hyperspaces of cones and fans}, Proc. Amer. Math. Soc, \textbf{77} (1979), no. 2, 279--288.

\smallskip

\bibitem{Esty} N. Esty, \textit{On the contractibility of certain hyperspaces}, Top. Proc., \textbf{32} (2008)


\bibitem{Illanes1}  A. Illanes, \textit{The hyperspace $C_{2}(X)$ for a finte graph is unique},
Glasnik Matematicki,  \textbf{37} (2002), 347--363.

\smallskip

\bibitem{Illanes2}  A. Illanes, \textit{Finite graphs $X$ have unique hyperspaces $C_{n}(X)$},
Top. Proc.,  \textbf{27} (2003), 179-188.

\smallskip

\bibitem{IllanesNadler} A. Illanes, S. Nadler, \textit{Hyperspaces: Fundamentals and Recent Advances}, Marcel Dekker, Inc., New York, 1999.

\end{thebibliography}

\end{document}